\documentclass[12pt]{amsart}
\usepackage{amscd}
\usepackage{verbatim}
\usepackage{amssymb, amsmath, amsthm, amscd,ifthen}
\usepackage[dvips]{graphics}
\usepackage[cp866]{inputenc}
\usepackage{graphicx}
\usepackage{epsfig}

\usepackage{amsmath,amssymb,amscd,}
\usepackage[all]{xy}

\newtheorem{theorem}{Theorem}

\newtheorem{lemma}{Lemma}
\newtheorem{remark}{Remark}
\newtheorem{dfn}{Definition}
\newtheorem{Ex}{Example}
\newtheorem{cor}{Corollary}

\title[Volume and lattice points  for  cyclopermutohedron]{Volume and lattice points counting for the cyclopermutohedron}
\author{Ilya Nekrasov, Gaiane Panina}

\begin{document}

\maketitle
\begin{abstract}

  The face lattice of the \textit{permutohedron}
realizes the combinatorics of linearly ordered partitions of the set
$[n]=\{1,...,n\}$. Similarly, the \textit{cyclopermutohedron} is a
virtual polytope that realizes the combinatorics of cyclically
ordered partitions of  $[n]$.

 It is known that the volume of the standard
permutohedron equals the number of trees with $n$ labeled vertices
multiplied by $\sqrt{n}$. The number of integer points of the
standard permutohedron equals the number of forests on $n$ labeled
vertices.

In the paper we prove that  the volume of the cyclopermutohedron
also  equals some weighted number of forests, which eventually
reduces to zero. We also derive a combinatorial formula for the
number of integer points in the cyclopermutohedron.

Another  object of the paper is the \textit{configuration space of
a polygonal linkage $L$}. It has a cell decomposition
$\mathcal{K}(L)$ related to the face lattice of cyclopermutohedron.
Using this relationship, we introduce and compute the volume
$Vol(\mathcal{K}(L))$.

\end{abstract}
\section{Introduction}\label{SectIntro}

The \textit{standard permutohedron} $\Pi_n$ is  defined (see
\cite{z}) as the convex hull of all points in $\mathbb{R}^n$ that
are obtained by permuting the coordinates of the point
$(1,2,...,n)$. It has the following properties:
\begin{enumerate}
    \item $\Pi_n$ is an $(n-1)$-dimensional polytope.
    \item The $k$-faces of  $\Pi_n$
are labeled by ordered partitions of the set
\newline $[n]=\{1,2,...,n\}$  into
$(n-k)$ non-empty parts.

    \item A face $F$ of $\Pi_n$ is contained in a face $F'$  iff the label of
$F$  refines the label  of $F'$.

Here and in the sequel, we mean the order-preserving refinement. For
instance, the label  $(\{1,3\},\{5,6\},\{4\},\{2\})$ refines the
label $(\{1,3\},\{5,6\},\{2,4\})$, but does not refine
$(\{1,3\},\{2,4\},\{5,6\})$.

    \item The permutohedron is a \textit{zonotope}, that is,  Minkowski sum
of line segments $q_{ij}$, whose defining vectors are
$\{e_i-e_j\}_{i<j}$, where $e_i$ are the standard basis vectors.

  \item The permutohedron splits into the union of \textit{bricks}
  (that is, some
elementary parallelepipeds)  labeled by all possible trees on $n$
vertices.
 The volume of each of the bricks equals
 $\frac{1}{\sqrt{n}}$, so $$Vol(\Pi_n)=\frac{1}{\sqrt{n}}\cdot
 \hbox{ number of trees on $n$ labeled vertices }=\sqrt{n}\cdot{n^{n-3}}{}.$$

 \item The number of integer points of the
standard permutohedron equals the number of forests on $n$ labeled
vertices. This fact  comes from some more delicate splitting of
$\Pi_n$ into
 bricks: unlike volume computing, we have to take
into account pieces of all dimensions, so we deal with
\textit{semiopen bricks} (details are given in Section
\ref{SectInteger}).
\end{enumerate}

 Similarly, the \textit{cyclopermutohedron} $\mathcal{CP}_{n+1}$ \cite{pan3} realizes the combinatorics of cyclically
ordered partitions of  $[n+1]=\{1,...,n, n+1\}$: all the $k$-faces
of the cyclopermutohedron are labeled by (all possible) cyclically
ordered partitions of the set $[n+1]$ into $(n+1-k)$ non-empty
parts, where $(n+1-k)>2$. The incidence relation in
${\mathcal{CP}}_{n+1}$ corresponds to the refinement: a cell $F$
contains a cell $F'$ whenever the  label of $F'$ refines
 the  label of $F$.

 The cyclopermutohedron is defined explicitly, as a weighted Minkowski sum of line segments.

In the paper we prove that  the volume of the cyclopermutohedron
  equals some weighted number of forests. Making use of the
theory of Abel polynomials, we eventually reduce the expression to
zero. We also give a combinatorial formula for the number of integer
points in the cyclopermutohedron.

Another  object of the paper is the \textit{configuration space, or
moduli spaces of a polygonal linkage $L$}. One of the motivations
for introducing the cyclopermutohedron is that  $\mathcal{CP}_{n+1}$
is a "universal" polytope for moduli spaces of polygonal linkages.
Namely, given a flexible polygon $L$, the space of its planar shapes
(that is, the configuration space) has a cell decomposition $\mathcal{K}(L)$,
whose combinatorics embeds in the combinatorics of   the face poset
of cyclopermutohedron. Using this relationship we introduce and
compute the volume $Vol(\mathcal{K}(L))$.

\bigskip

The paper is organized as follows. In Section \ref{SecTheorBackgr}
we give all necessary information about virtual polytopes, and also
the definition and properties of the cyclopermutohedron. Abel
polynomials are also sketched in the section.

In Section \ref{SecVolIsZero} we explain the meaning of the ''volume
of cyclopermutohedron'', and prove that it
 equals zero.  In Section
\ref{section_linkages} we explain the relationship with polygonal
linkages and give a formula for the volume of the configuration
space (Theorem \ref{ThmVolLink}).

 Finally, in Section
\ref{SectInteger}  we compute the number of integer points in the
cyclopermutohedron (Theorem \ref{ThmInteger}).

\bigskip

\textbf{Acknowledgements.} The present research  is  supported by
RFBR, research project No. 15-01-02021. The first author was also
supported
 by the Chebyshev Laboratory under RF Government grant
11.G34.31.0026, and  JSC ''Gazprom Neft''.

\section{Theoretical backgrounds}\label{SecTheorBackgr}
\subsection{Virtual polytopes}
Virtual polytopes appeared in the literature as useful
geometrization of Minkowski differences of convex polytopes. A
detailed discussion can be found in \cite{pkh,pan2,panstr}; below we
give just a brief sketch. As a matter of fact, in the paper (except
for Section \ref{section_linkages}) we need no geometrization. Even
for volume and integer point counting, it is sufficient to know that
virtual polytopes form the Grothendieck group associated to the
semigroup  of convex polytopes.

More precisely, \textit{a convex polytope }is the convex hull of a
finite, non-empty point set in the Euclidean space $\mathbb{R}^n$.
Degenerate polytopes are also included, so a closed segment and a
point are polytopes, but not the empty set. We denote by
$\mathcal{P}^+$ the set of all convex polytopes.

Let $K$ and $L \in \mathcal{P}^+$ be  two convex polytopes.  Their
\textit{Minkowski sum} $K + L$ is defined by:
$$
K + L = \{\textbf{x} +\textbf{y} : \textbf{x} \in K, \textbf{y} \in
L\}.
$$

Minkowski addition turnes  the set $\mathcal{P}^+$ to  a commutative
semigroup whose unit element is the convex set containing exactly
one point \newline $E= \{ 0 \}$.

\begin{dfn}
The group $\mathcal{P}$ of {\em virtual polytopes} is the
Grothendick group associated to the semigroup $\mathcal{P}^+$ of
convex polytopes under Minkowski addition.

The elements of $\mathcal{P}$ are called {\em virtual polytopes}.
\end{dfn}
More instructively,  $\mathcal{P}$ can be explained as follows.

\begin{enumerate}
  \item A virtual polytope is
 a formal difference $K- L$.

  \item Two such
expressions $K_1- L_1$ and $K_2- L_2$ are identified, whenever $K_1+
L_2=K_2+ L_1$.
  \item  The group operation is
 defined by  $$(K_1- L_1) + (K_2- L_2):=
(K_1 + K_2)- (L_1 + L_2).$$

\end{enumerate}

It is important that the notions of ''volume'' and  ''number of
integer points'' extend nicely  to virtual polytope. We explain
these constructions in the subsequent sections.
\subsection{Cyclopermutohedron}\cite{pan3}

Assuming that $\{e_i\}_{i=1}^n$ are standard basic vectors in
$\mathbb{R}^n$, define the points

$$\begin{array}{ccccccccc}
    R_i=\sum_{i=1}^n (e_j-e_i)=(-1, & ... & -1, & n-1, & -1, & ... & -1, & -1, &-1, )\in \mathbb{R}^{n},\\
     &  & & \ i &  &  &  &  &
  \end{array}
$$
and  the following two families of line segments:
$$q_{ij}=\left[e_i,e_j\right], \ \ \  i<j$$ and $$  r_i=\left[0,R_{i} \right].$$

We also need the point $e=\left(1,1,...,1\right)\in \mathbb{R}^{n}$.

 The  \textit{cyclopermutohedron} is a virtual polytope defined as    the  Minkowski sum:
$$ \mathcal{CP}_{n+1}:=   \bigoplus_{i< j} q_{ij} + e- \bigoplus_{i=1}^n  r_i.$$

Here and in the sequel, the sign ''$ \bigoplus$'' denotes the
Minkowski sum, whereas  the sign ''$\sum$'' is reserved for the sum
of numbers.

\bigskip

The cyclopermutohedron $ \mathcal{CP}_{n+1}$
 lies in the hyperplane
 $$x_1+...+x_n=\frac{n(n+1)}{2},$$
 so its actual dimension is $(n-1)$.

\begin{remark}\label{RemPermSum} The Minkowski sum
$$   \bigoplus_{i< j} q_{ij}+ e$$  is known to be equal to the standard permutohedron
$\Pi_n$  (see \cite{z}). Therefore we can write

$$ \mathcal{CP}_{n+1}=  \Pi_n - \bigoplus_{i=1}^n  r_i.$$
\end{remark}

The face poset of $\mathcal{CP}_{n+1}$ encodes cyclically ordered
partitions of  the set $[n+1]=\{1,...,n+1\}$:
 \begin{enumerate}
   \item For  $k=0,...,n-2$, the   $k$-dimensional faces  of ${\mathcal{CP}}_{n+1}$ are labeled by (all possible)  cyclically ordered
partitions of the set $[n+1]$  into $(n-k+1)$ non-empty parts.
   \item A face $F'$ is a face of $F$   whenever the  label of $F'$ refines
 the  label of $F$. Here  we mean order
 preserving refinement.
 \end{enumerate}

\subsection{Abel polynomial and rooted forests}\cite{Sagan}\label{SecAbel}

A \textit{rooted forest } is a graph equal to a  disjoint union of
trees, where each of the trees has a marked vertex.

The\textit{ Abel polynomials  } form a  sequence of polynomials,
where the $n$-th term is defined by

$$A_{n,a}(x)=x(x-an)^{n-1}.$$

 A
special case of the Abel polynomials with $a=-1$ counts rooted
labeled forests. Namely, if  $A_{n}(x) := A_{n,-1}(x) =
 x(x+n)^{n-1}$ is the $n$-th {Abel polynomial}, then

$$A_{n}(x)=\sum_{k = 0}^{n} t_{n,k}\cdot x^{k}  ,$$
  where $ t_{n,k}$
is the number of forests on $n$ labeled vertices consisting of $k$
rooted trees.

\section{Volume of cyclopermutohedron equals zero}\label{SecVolIsZero}

The notion of volume extends nicely from convex polytopes to virtual
polytopes. We explain below the meaning of the   \textit{volume of a
virtual zonotope}.

 Assume
we have a convex zonotope $Z\subset \mathbb{R}^n$, that is, the
Minkowski sum of some linear segments $\{s_i\}_{i=1}^m$:
$$Z=\bigoplus_{i=1}^m \  s_i.$$

For each subset $I\subset [m]$ such that $|I|=n$, denote by $Z_I$
the \textit{elementary parallelepiped}, or the \textit{brick}
spanned by $n$ segments $\{s_i\}_{i\in I}$, provided that the
defining vectors of the segments are linearly independent.  In other
words, the brick equals the Minkowski sum $$Z_I=\bigoplus_Is_i.$$

It is known that $Z$ can be partitioned into the union of all such
$Z_I$, which implies immediately
$$Vol(Z)=\sum_{I\subset [m], |I|=n}Vol(Z_I)=\sum_{I\subset [m], |I|=n}|Det(S_I)|,$$
 where $S_I$ is the matrix composed of defining vectors of the
 segments from $I$.

 Now take  positive $\lambda_1,...,\lambda_m$ and  sum up the dilated segments $\lambda_is_i$. Clearly, we have

 $$Vol\Big(\bigoplus_{i=1}^m\ \lambda_i s_i\Big)=\sum_{I\subset [m], |I|=n}\prod_{i\in I}\lambda_i \cdot|Det(S_I)|.$$

 For fixed $s_i$, we get a polynomial in
 $\lambda_i$, which counts  not only the volume of convex zonotope (which originates from positive $\lambda_i$),
 but also the volume of a virtual zonotope, which originates from any real $\lambda_i$, including negative ones, see \cite{pkh,panstr}.

 So, one can use the above formula as the definition of the volume
 of a virtual zonotope.

An almost immediate consequence is:
 \begin{lemma}\label{lemmaVol}Let $E=E_n$ be the set  of edges of the complete graph $K_n$.
The \newline $(n-1)$-volume of the cyclopermutohedron can be
computed by the formula:
   $$Vol(\mathcal{CP}_{n+1})=  Vol\Big( \bigoplus_{i< j} q_{ij} - \bigoplus_{i=1}^n
 r_i\Big)=$$
 $$=\frac{1}{\sqrt{n}}\sum_{|I|+|M|=n-1} (-1)^{|M|}|Det(q_{ij},r_k,e)|_{(ij)\in I, \ k \in M}.$$
 Here $I$ ranges over  subsets of $E$, whereas $M$ ranges over  subsets of
 $[n]$. The matrix under determinant is composed of defining vectors of the segments
 $q_{ij}$ and
 $r_k$, and also of the  vector $e=(1,1,...,1,1)$.
 \end{lemma}

\textit{ Proof.} The cyclopermutohedron $\mathcal{CP}_{n+1}$  lies in the
hyperplane

$$x_{1} + \dots + x_{n} = \frac{n(n+1)}{2},$$

so its dimension equals  $ n - 1$. That is, we deal with $(n-1)$-volume, which reduces to the $n$-volume by adding $e=(1,1,...,1,1,)$
and dividing by   $|e|=\sqrt{n}$.\qed

\begin{remark} \label{RemBrGram} The formula for the volume of a virtual zonotope also has a geometrical meaning which we
 briefly sketch here.   Due to Brianchon-Gram decomposition of virtual polytopes (see \cite{panstr} or \cite{pkh}), any virtual polytope can be viewed as a codimension one homological
 cycle, and therefore possesses a well-defined (algebraic) volume.

For a virtual zonotope, the associated cycle  decomposes into
homological sum of elementary bricks, but the latter should be
understood also as homological cycles coming  with different
orientations.
 More precisely, if the number
of negative $\lambda_i$ in the sum $\bigoplus_{i=1}^n \lambda_i s_i$
is even, then the corresponding elementary brick  equals the
boundary of elementary parallelepiped   $\ \ \partial
\Big(\bigoplus_{i=1}^n |\lambda_i| s_i\Big)$  with the
\textit{positive} orientation (that is, cooriented by the outer
normal vector). If the number of negative $\lambda_i$ is odd, we
have the same cycle with the \textit{negative} orientation.
\end{remark}

\begin{theorem}\label{TeoremVol}
$Vol(\mathcal{CP}_{n+1}) = 0$.
\end{theorem}

Proof. Keeping in mind  Lemma \ref{lemmaVol}, let's first fix $I$
and $M$ with $|I|+|M|=n-1$, and compute one single  summand
$|Det(q_{ij},  r_{k},  e)|_{(ij)\in I, \ k \in M}$.

If $M=\emptyset$, the determinant equals $1$ iff the set $I$ gives a
tree. Otherwise it is zero. (This is the reason for  the volume
formula of the permutohedron.)

Assume now that $M$ is not empty.

$$|Det(q_{ij},  r_{k},  e)| =
\begin{vmatrix} 0 & 0 & \dotsm & -1 & \dotsm & 1\\
                \vdots & \vdots & \ddots & -1 & \dotsm & 1\\
                -1  & 0 & \dotsm & -1 & \dotsm & 1\\
                \vdots & -1 & \dotsm & -1 & \dotsm & 1\\
                1   & 0 & \dotsm & -1 & \dotsm & 1\\
                \vdots & \vdots & \ddots & -1 & \dotsm & 1\\
                0 & 0 & \dotsm & n-1 & \dotsm & 1\\
                \vdots & \vdots & \ddots & -1 & \dotsm & 1\\
                0   & 1 &\dotsm & -1 & \dotsm & 1\\
                \vdots & \vdots & \ddots & -1 & \dotsm & 1\\
                0 & 0 & \dotsm & -1 & \dotsm & 1\\
\end{vmatrix} = $$
Adding  $e$ to all the columns $r_i$, we get:

$$ = n^{|M|} \cdot \begin{vmatrix} 0 & 0 & \dotsm & 0 & \dotsm & 1\\
                \vdots & \vdots & \ddots & 0 & \dotsm & 1\\
                -1  & 0 & \dotsm & 0 & \dotsm & 1\\
                \vdots & -1 & \dotsm & 0 & \dotsm & 1\\
                1   & 0 & \dotsm & 0 & \dotsm & 1\\
                \vdots & \vdots & \ddots & 0 & \dotsm & 1\\
                0 & 0 & \dotsm & 1 & \dotsm & 1\\
                \vdots & \vdots & \ddots & 0 & \dotsm & 1\\
                0   & 1 &\dotsm & 0 & \dotsm & 1\\
                \vdots & \vdots & \ddots & 0 & \dotsm & 1\\
                0 & 0 & \dotsm & 0 & \dotsm & 1\\
\end{vmatrix} =n^{|M|}\cdot(*).$$
We wish to proceed in a similar way, that is, add the columns
containing the unique entry $1$ to other columns chosen in an
appropriate way.
 To explain this reduction
 let us give two technical definition.

\bigskip
\begin{dfn}\label{DefDecForest}

A \textit{decorated forest} $F=(G,M)$  is a graph $G=([n],I)$
without cycles on $n$ labeled vertices together with a set of marked
vertices $M\subset[n]$ such that the following conditions
hold:\begin{enumerate}
    \item Number of marked vertices $|M|$  $+$ number of edges  $|I|$
    equals $n-1$.
    \item Each connected component of $G$ has at most one marked
    vertex.
\end{enumerate}
\end{dfn}

Immediate observations are:

 Each decorated forest has exactly one
connected component with no vertices marked. We call it\textit{ a
free tree}. Denote by $N(F)$ the number of vertices of the free
tree.

 Each decorated forest is a disjoint union of the free tree and some
rooted forest. The number of rooted trees equals $|M|$.

 Each decorated forest $F$  yields a
collection of $\{e_{ij}, r_k\}_{(ij)\in I, \ k\in M}$, whose above
determinant $(*)$  we denote by $|Det(F)|$  for short.

\begin{figure}[h]
\centering
\includegraphics[width=12 cm]{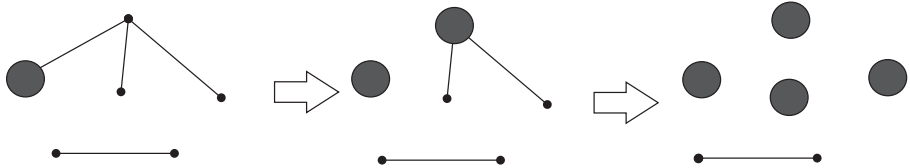}
\caption{Reduction process for a forest with $N(F)=2,\ |M|(F)= 1$.
Grey balls denote the marked vertices.}\label{FigKill}
\end{figure}

For instance, for the first decorated forest in Figure
\ref{FigKill}, we have  $N(F)=2,\ |M|= 1.$

\bigskip

 Now we define the
\textit{reduction of a decorated forest} (see Figure \ref{FigKill}
for example). It goes as follows.

Assume we have a decorated forest.  Take a marked vertex $i$ and
 an incident edge $(ij)$. Remove the edge and mark the vertex $j$.
Repeat until is possible.

Roughly speaking, a marked vertex $i$ {kills} the edge $(ij)$ and
{generates} a new marked vertex $j$.
%\end{dfn}

\bigskip

An obvious observation is:
\begin{lemma}\begin{enumerate}
    \item  The free tree does not change during the reduction.
    \item The reduction brings us to a
decorated forest with a unique  free tree.  All other trees are
one-vertex trees, and all these vertices are marked.
    \item The reduction can be shortened: take the
    connected components one by one and do the following.\begin{enumerate}
        \item If a connected component has no marked vertices, leave
        it as it is.
        \item If a connected component has a marked vertex,
        eliminate all its edges and mark all its vertices.
    \end{enumerate}
    \item The reduction does not depend on the order of the marked
 vertices we deal with.\qed
\end{enumerate}

\end{lemma}

Before we proceed with the proof of Theorem \ref{TeoremVol}, prove
the lemma:
\begin{lemma}\label{LemmaSmallDet}\begin{enumerate}
    \item For each decorated forest $F$, $$|Det(F)|= N(F).$$

    \item If a collection $\{e_{ij},
r_k\}$ does not come from a decorated forest, that is, violates
condition (2) from Definition \ref{DefDecForest}, then
$$|Det(e_{ij},
r_k)|= 0.$$
\end{enumerate}
\end{lemma}

\textit{Proof of the lemma.} (1) For a decorated forest, we
manipulate with the columns according to the reduction process. We
arrive at a matrix which (up to a permutation of the columns and up
to a sign) is:

$$\left(
    \begin{array}{ccc}
      A & O &1\\
      O & E &1\\
    \end{array}
  \right).
$$

Here $A$ is the matrix corresponding to the free tree, $E$ is the
unit matrix, and the very last column is $e$. Its determinant equals
$1$.

(2) If the collection of vectors does not yield a decorated forest,
that is, there are two marked vertices on one connected component,
the analogous reduction gives a zero column.\qed

\bigskip

 Basing on Lemmata
\ref{LemmaSmallDet} and \ref{lemmaVol}, we  conclude:

$$Vol( \mathcal{CP}_{n+1}) = \frac{1}{\sqrt{n}}\,\sum_{F } (-n)^{|M(F)|} \cdot N(F)=, $$

where the sum extends over all decorated forests $F$ on $n$
vertices. (Remind that  $M(F)$ is the set of marked vertices,
 $N(F)$ is the number of
vertices of the free tree.)

Next, we group the forests by the number $N=N(F)$ and write

$$= \frac{1}{\sqrt{n}}\sum_{ N = 1}^{n} \binom{n}{N} N^{N-2} \cdot N \sum_f (-n)^{C(f)} =$$
$$= \frac{1}{\sqrt{n}}\sum_{ N = 1}^{n} \binom{n}{N} N^{N-1}\sum_f (-n)^{C(f)}=\frac{1}{\sqrt{n}}\cdot(**),$$

where the second sum ranges over all rooted forest on $(n-N)$
labeled vertices, $C(\cdot)$ is the number of connected components.

Let us explain this in more details.\begin{enumerate}
    \item $N$ ranges from $1$ to $n$. We choose $N$ vertices
  in  $\binom{n}{N}$ different ways and place a tree on these
  vertices in $ N^{N-2}$ ways.
    \item On the rest of the vertices we  place a rooted forest.
\end{enumerate}

Recalling that   $t_{n-N,k}$ the number of forests on $(n-N)$
labeled vertices of $k$ rooted trees, we write:
$$(**)= \sum_{ N = 1}^{n} \binom{n}{N} N^{N-1}\sum_{k=1}^{n-N} (-n)^{k}\cdot t_{n-N,k}.$$

 Section \ref{SecAbel} gives us:

$$\sum_{k = 0}^{n} t_{n,k}x^{k} = A_{n}(x),$$

where $A_{n}(x) = x(x+n)^{n-1}$ is the {Abel
polynomial}.

Setting  $-n = x$, we get

$$\sum_{k=1}^{n-N} (-n)^{k}\cdot t_{n-N,k} =  A_{n-N}(-n).$$

Thus $(**)$  converts to

%$$\sum_{N = 1}^{n}\binom{n}{N} N^{N - 1}A_{n - N}(-1,-n) =: Q_{n}\textbf{????????}$$

%or

$$\sum_{N = 1}^{n}\binom{n}{N} N^{N - 1}A_{n - N}(-n) =: Q_{n}.$$

Applying the definition of $A_{n-N}(-n)$, we get

$$Q_{n}  = \sum_{N = 1}^{n}\binom{n}{N} N^{N - 1} (-n)(- n + n - N)^{n-N-1} =  (-1)^{n} n \cdot \sum_{N = 1}^{n}(-1)^{N}\binom{n}{N} N^{n - 2}  .$$

Introduce the following polynomial:
$$p(x) := \sum_{N = 0}^{n}N^{n - 2}\binom{n}{N}  x^{N},$$
for which we have $Q_{n} = p(-1)$. Set also

$$p_{0}(x) := (1 + x)^{n} = \sum_{N = 0}^{n} \binom{n}{N} x^{N},$$
$$p_{i}(x):= x\cdot p'_{i-1}(x) = \sum_{N = 0}^{n} N^i\binom{n}{N}
 x^{N}.$$ We clearly have $p(x) = p_{n-2}(x)$.
Besides,  $(1+x)^{n - k}$ divides $p_{k}(x)$,  which implies $Q_{n}
= 0$. \qed

\section{Polygonal linkages: volume of the configuration space}\label{section_linkages}
\subsection{Definitions and notation}\label{subsection_linkages_notation} A \textit{flexible $(n+1)$-polygon}, or a \textit{polygonal  $(n+1)$-linkage} is
a sequence of positive numbers $L=(l_1,\dots ,l_{n+1})$. It should
be interpreted as a collection of rigid bars of lengths $l_i$ joined
consecutively in a closed chain by revolving joints. We always
assume that the triangle inequality holds, that is, $$\forall j, \ \
\ l_j< \frac{1}{2}\sum_{i=1}^{n+1} l_i$$ which guarantees that the
chain of bars can close.

We also assume that the last bar is the longest one:
$$\forall j \ \ \ l_{n+1} \geq l_j.$$

 \textit{A planar configuration} of $L$   is a sequence of points
$$P=(p_1,\dots,p_{n+1}), \ p_i \in \mathbb{R}^2$$ with
$l_i=|p_ip_{i+1}|$, and $l_{n+1}=|p_{n+1}p_{1}|$. As follows from
the definition, a configuration may have self-intersections and/or
self-overlappings.

 \textit{The moduli space, or the
configuration space $M(L)$}  is the space  of all configurations
modulo orientation preserving isometries of $\mathbb{R}^2$.

Equivalently, we can define $M(L)$ as
$$M(L)=\{(u_1,...,u_{n+1}) \in (S^1)^{n+1} : \sum_{i=1}^{n+1} l_iu_i=0\}/SO(3).$$

The (second) definition shows that $M(L)$ does not depend on the
ordering of $\{l_1,...,l_{n+1}\}$; however, it does depend on the
values of $l_i$.

Let us comment on this dependance. Consider $(l_1,...,l_{n+1})$ as a
point in the parameter space $\mathbb{R}^{n+1}$. The hyperplanes in
$\mathbb{R}^{n+1}$ defined by all possible equations
$$\sum_{i=1}^{n+1}\varepsilon_il_i=0 \hbox{ \ \ with \ \ } \varepsilon_i=\pm
1$$ are called\textit{ walls}. Throughout the section we assume that
 the
point $\{l_1,...,l_{n+1}\}$ belongs to none of the walls. This
genericity assumption implies that the moduli space $M(L)$ is a
closed $(n-2)$-dimensional manifold.

 The walls  dissect
$\mathbb{R}^{n+1}$ into a number of \textit{chambers}; the topology
of $M(l_1,...,l_{n+1})$ depends only on the chamber containing
$\{l_1,...,l_{n+1}\}$ (see \cite{F}).

The  manifold $M(L)$ is already well studied. In this paper we make
use of  the described below cell structure on the space $M(L)$.

\subsection{ The complex $\mathcal{K}(L)$}\label{SectionCW}

Assume that $(l_1,...,l_{n+1})$ is fixed.

 A set $I\subset [n+1]=\{1,2,...,n+1\}$ is called
\textit{short}, if
$$\sum_{i \in I}^{}l_i <\frac{1}{2} \sum_{i=1}^{n+1}l_i.$$
Otherwise $I$ is a\textit{ long} set.

A partition of the set  $[n+1]$ is called \textit{admissible} if
all the sets in the partition are short.

\begin{theorem}\cite{pan2}\label{ThmCellComplex}
There is  a structure of a regular CW-complex $\mathcal{K}(L)$ on
the moduli space $M(L)$. Its complete combinatorial description
reads as follows:
\begin{enumerate}
    \item  $k$-cells of the complex $\mathcal{K}(L)$ are labeled by cyclically ordered admissible partitions of
 the set  $[n+1]$  into $(n+1-k)$ non-empty
parts.

    \item A closed cell $C$ belongs to the boundary of some other closed cell
    $C'$  iff  the partition  $\lambda(C)$ is finer than
    $\lambda(C')$.

\end{enumerate}
\end{theorem}

{\bf A remark on notation.} We write a cyclically ordered partition
as a (linearly ordered) string of sets where the set containing the
entry ''$n$'' stands on the last position.

We stress  that the order of the sets matters, whereas there is no
ordering inside a set. For example,
$$(\{1\} \{3 \} \{4,  2, 5,6\})\neq(\{3 \}\{1\}  \{4,  2, 5,6\})= ( \{3 \}\{1\}\{ 2,4, 5,6\}).$$

\begin{Ex}\label{ExPErmuto}
Assume that
$$l_{n+1}=\sum_{i=1}^{n}l_i-\varepsilon,$$
where $\varepsilon$ is small. In this case the moduli space $M(L)$
is the sphere $S^{(n-2)}$, see \cite{F}, and the complex
$\mathcal{K}(L)$ is isomorphic to the boundary complex of the
permutohedron $\Pi_{n}$.
\end{Ex}

For any $(n+1)$-linkage $L$, the complex
 $\mathcal{K}(L)$ automatically embeds  in the face complex of cyclopemutohedron $\mathcal{CP}_{n+1}$, and therefore can be realized by a polyhedron
 which we denote by $\mathcal{P}(L)$.   Vividly speaking, the polyhedron $\mathcal{P}(L)$ is patched of those faces of the cyclopermutohedron, whose labels are admissible partitions.

\medskip
\textbf{Example.}

 For $L$ as in Example \ref{ExPErmuto},
  $\mathcal{P}(L)$ equals the boundary  of the permutohedron $\Pi_{n}$.

\medskip
\textbf{Example.}

 Let $n=5$, $L=(1.2,\ 1,\ 1,\ 0.8,\ 2.2)$. Then $\mathcal{P}(L)$ is the cylinder depicted in
Fig. \ref{torus}. The two shadowed faces are labeled by
$(\{123\}\{4\})$  and $(\{4\}\{123\})$. Since the partitions
$(\{123\}\{4\}\{5\})$  and $(\{4\}\{123\}\{5\})$ are non-admissible,
 these faces of permutohedron are removed, whereas all other faces
of the permutohedron survive. There are also six  "diagonal"
rectangular faces. They are labeled by $(\{12\}\{3\}\{45\})$,
$(\{1,3\}\{2\}\{45\})$, $(\{2\}\{13\}\{45\})$,
$(\{23\}\{1\}\{45\})$, $(\{3\}\{12\}\{45\})$, and
$(\{1\}\{23\}\{45\})$.

\begin{figure}[h]\label{torus}
\centering
\includegraphics[width=12 cm]{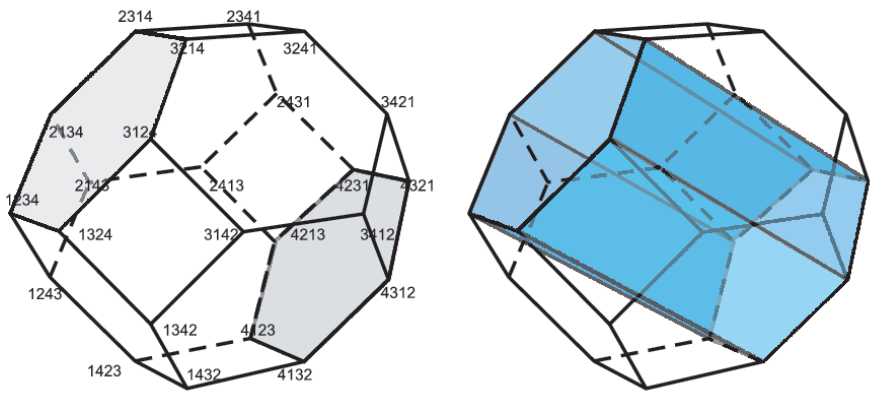}
\caption{The complex $\mathcal{K}(L)$ for the $5$-linkage \newline
$L=(1,2;\ 1;\ 1;\ 0,8;\ 2,2)$. We remove from the permutohedron the
two shaded facets and patch in the  cylinder.}\label{torus}
\end{figure}

\subsection{Volume of the complex $\mathcal{K}(L)$} Following the ideology of Remark \ref{RemBrGram}, $\mathcal{P}(L)$ can be viewed as a codimension one homological cycle
(or as a generalization of  closed piecewise linear oriented
manifold) in the Euclidean space. Therefore it makes sense to speak
of the \textit{volume} of the part of the space bounded  by
$\mathcal{P}(L)$.
 Since $\mathcal{P}(L)$  may have  many self-intersections,  the volume  means the \textit{algebraic volume}, that is,
multiplicities  (which can be also negative) are taken into account.

Let us explain this in more details. For each point $x\in
\mathbb{R}^{n}$, denote by $ind_x(\mathcal{P}(L))$ the index of the
cycle with respect to the point $x$. Then by the  \textit{volume of
the configuration space} we mean
$$Vol(M(L)):=Vol(\mathcal{P}(L)):=\int_{\mathbb{R}^n}ind_x(\mathcal{P}(L))dx.$$

\begin{dfn}For an $(n+1)$-linkage $L$, a  decorated forest $F$ on $n$ labeled vertices is called \textit{non-admissible}, if  the vertex set of  the free
 tree is a long set.
\end{dfn}
In notation of Section \ref{SecVolIsZero}, the following lemma
holds:
\begin{lemma}\label{lemmaVolLink} For a $(n+1)$-linkage $L$, we have: $$Vol(M(L) ) = \frac{1}{\sqrt{n}}\,\sum_{non-admissible \ F } (-n)^{|M(F)|} \cdot N(F), $$
where the sum ranges over all non-admissible decorated  forests on
$n$ labeled vertices.

We remind that $|M(F)|$ denotes the number of marked vertices,
$N(F)$ is the number of vertices of the free tree.
\end{lemma}

Proof.
 Let us take the linkage
$L_0=(l_1,...,l_{n},\lambda)$ assuming that the value of $\lambda$
continuously and monotonly changes from
$$\sum_{i=1}^{n}l_i-\varepsilon \hbox{ \ \ \  to \ \ \  } l_{n+1}.$$

In the beginning   we have the permutohedron $\Pi_n$, whose volume
we already know. At the end, we have $\mathcal{P}(L)$, whose volume
we wish to calculate. In between we have a (finite) number of Morse
surgeries, and we can control the behavior of the volume at each of
the surgeries.

Prove first that the formula holds true for
$\lambda=\sum_{i=1}^{n}l_i-\varepsilon$. Indeed, for this particular
$\lambda$, a "decorated non-admissible forest on $n$ vertices" means
just "a free tree on $n$ vertices", so the statement of the theorem
reduces to the formula for the volume of the standard permutohedron,
see Section \ref{SectIntro}.

Now we start changing $\lambda$. This means that we have a path in
the parameter space $\mathbb{R}^{n+1}$, which  crosses some of the
walls. We can assume that the walls are crossed one by one; if this
is not the case, we perturb generically the original lengths $l_i$.

 Once we
cross a wall, the complex $\mathcal{K}$, and its polytopal
realization change by a  surgery which we describe below. Denote by
$Pol_{Old}$ and by $Pol_{New}$ the polyhedra that realize
$\mathcal{K}$ before and after the surgery respectfully.

Let us look at the surgery in more details. Once a  wall is crossed,
 some maximal by
inclusion short set $T\subset \{1,...,n\}$ turns to a long set,
whereas its complement $\overline{T}=[n+1]\setminus T$ becomes
short. We conclude that the new complex $\mathcal{K}$ can be
obtained from the old complex by removing some of the cells and
adding some new cells. The cells that get removed after crossing the
wall are labeled by $(*,T,*)$. Here
 whereas the new cells that
appear are labeled by $(*,[n]\setminus T)$. Here ''$*$'' means just
any ordered partition of the complement assuming that altogether we
have at least three parts.

 The cells that get removed form a subcomplex isomorphic to  the boundary of the permutohedron  $\Pi_{n-|T|}$
    multiplied by a $(|T|-1)$-ball.
     The cell structure of $\mathcal{K}$ converts this ball to the permutohedron $\Pi_{|T|}$.
    So, we have the following Morse surgery: we cut out the cell subcomplex $$\mathcal{C}_1=(\partial \Pi_{n-|T|})\times \Pi_{|T|},$$ and  patch instead the cell complex $$\mathcal{C}_2=\Pi_{n-|T|}\times \partial \Pi_{|T|}$$ along the identity mapping on
     their common boundary $\partial \Pi_{n-|T|}\times \partial \Pi_{|T|}$.

Denote by $\mathcal{C}:=\mathcal{C}_1 \cup \mathcal{C}_2$ the union
of these complexes. Combinatorially, we have $\mathcal{C}=\partial
\Big(\Pi_{n-|T|}\times
 \Pi_{|T|}\Big).$

  $\mathcal{C}$ (taken with an appropriate orientation) relates the old and new
polyhedra. Namely, we have a homological sum:

$$Pol_{New} =Pol_{Old} +\mathcal{C}.$$

 This means that the new and old volumes
are related by
$$Vol(Pol_{New}) =Vol(Pol_{Old}) +Vol(\mathcal{C}).$$

After geometrically realizing these complexes, we decompose the
realization of
 $\Pi_{n-|T|}\times
 \Pi_{|T|}$ into the homological sum of bricks $P_i\times P_j$, where  $P_i$ is an elementary brick from $\Pi_{n-|T|}$, and $P_j$ is an elementary brick from $\Pi_{|T|}$.
 The first elementary
brick $P_i$ corresponds to a tree on $T$, whereas $P_j$ corresponds
to a  tree on  $[n+1]\setminus T$, or, equivalently, to a rooted
forest  on $[n]\setminus T$. In other words, each such pair
$(P_i,P_j)$ gives us a rooted forest $F$  whose free tree is
non-admissible.

The brick $P_i\times P_j$ has a geometrical realization as the
Minkowski sum of corresponding line segments. It contributes
 $(-n)^{|M(F)|}$ to $Vol(\mathcal{C})$.

Therefore, if the statement of the theorem is true for $Pol_{Old}$,
it is also true for $Pol_{New}$.\qed

\begin{theorem}\label{ThmVolLink} For a flexible $(n+1)$-polygon $L$, we have: $$Vol(M(L) )={\sqrt{n}}\,\sum_{k=0}^n (-1)^{k} \cdot a_{k}\cdot (n-k)^{n-2},$$
where $a_k$   is the number of $(k+1)$-element short subsets of
$[n+1]$ containing the entry $(n+1)$.

\end{theorem}

Proof.  Using Lemma \ref{lemmaVolLink}, we first fix a number $k$
and choose a long $k$-element subset of $[n]$. This can be done in
$a_{n-k}$ ways. We put a tree on these vertices in $k^{k-2}$ ways
and arrive at

$$Vol(M(L) ) =
\frac{1}{\sqrt{n}}\,\sum_{k=1}^n a_{n-k}\cdot k^{k-2}\sum_{g \ is \
a \ rooted \ forest\  on\  (n-k)\  vertices} (-n)^{C(g)} \cdot
N(F)=$$$$=\frac{1}{\sqrt{n}}\,\sum_{k=1}^n a_{n-k}\cdot
k^{k-1}\sum_{g \ is \ a \ rooted \ forest\  on\  (n-k)\ vertices}
(-n)^{C(g)}=
$$

By the identity from Section \ref{SecAbel}
$$\sum_{g \ is \ a \ rooted \ forest\    on \ m \
vertices}x^{C(g)}=x\cdot(x+m)^{m-1},$$ we get
$$=\frac{1}{\sqrt{n}}\,\sum_{k=1}^n a_{n-k}(-k)^{k-1}\cdot(-n)\cdot
k^{n-k-1}
$$

$$=-\sqrt{n}\,\sum_{k=1}^n a_{n-k}k^{k-1}\cdot (-k)^{n-k-1}=
$$
$$=\sqrt{n}\,\sum_{k=1}^n a_{n-k}k^{n-2}\cdot (-1)^{n-k}.
$$ Interchanging $k$ and $n-k$, we get the desired.\qed

\bigskip

\textbf{Remark.}  Betty numbers  $\beta_k=\beta_k(M(L))$ are
expressed in terms of $a_k$, see \cite{faS}:
$$\beta_k=a_k+a_{n-k-3}.$$

\bigskip

\begin{cor}
 Assume $n+1=2m+1$. For  the equilateral
$(n+1)$-linkage $L=(1,1,...,1)$ we have:

$$Vol(M(L) ) =
{\sqrt{2m}}\,\sum_{k=0}^m (-1)^{k} \cdot \binom{2m}{k}\cdot
(2m-k)^{2m-2}.$$
\end{cor}

Proof. Indeed, for  the equilateral linkage, ''a short set'' means
''a set with cardinality $\leq m$''. Therefore
$$a_k=\left\{
                              \begin{array}{ll}
                                \binom{n-1}{k}, & \hbox{  if  $k\leq m$;} \\
                                0, & \hbox{ otherwise.\qed}
                              \end{array}
                            \right.$$

\section{Integer points counting for
cyclopermutohedron}\label{SectInteger}
\subsection{Integer points counting for cyclopermutohedron: theoretical backgrounds}
\begin{figure}[h]
\centering
\includegraphics[width=10 cm]{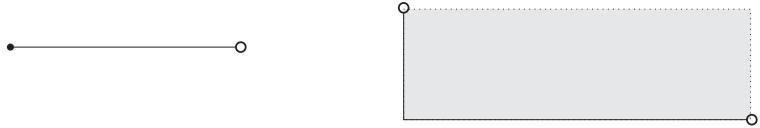}
\caption{A semiopen segment and a semiopen rectangle. The dashed
lines and white points are missing.}\label{FigSemiopen}
\end{figure}

The first leading idea for integer points enumeration in a zonotope
is to decompose it into elementary bricks, as we did in Section
\ref{SectIntro}.  However, unlike volume computation, we have to
take into account the ''pieces'' of all dimensions, including
points. By this reason we introduce \textit{semiopen bricks}. The
latter are Minkowski sums of semiopen segments, see Figure
\ref{FigSemiopen}.

\begin{figure}[h]
\centering
\includegraphics[width=10 cm]{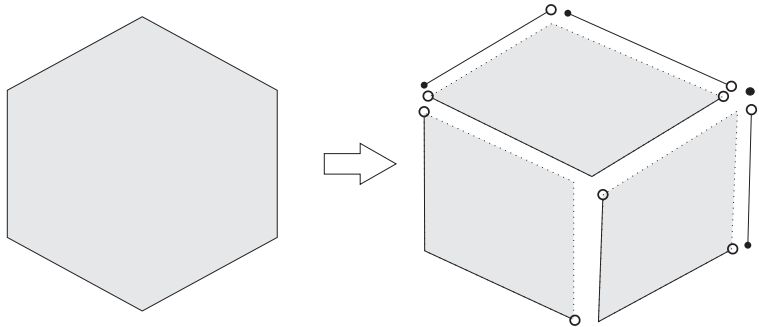}
\caption{ The permutohedron $\Pi_3$ splits into three semiopen
parallelograms, three semiopen segments, and one
point.}\label{FigSplit}
\end{figure}

A zonotope decomposes in a disjoint union of semiopen bricks of
dimensions ranging from $0$ to $n$.

\begin{Ex}Permutohedron $\Pi_n$ decomposes in a disjoint union of semiopen
bricks that are in a one-to-one correspondence with forests on $n$
labeled vertices. Each of the bricks contributes exactly one integer
point, so for the number of integer points $\Lambda$, we have:
$$\Lambda(\Pi_n)=\hbox{ number of forests
on $n$ labeled vertices.}$$
\end{Ex}

Below we almost literally repeat the arguments from Section
\ref{SecVolIsZero}. Assume we have a convex zonotope $Z\subset
\mathbb{R}^n$, that is, the Minkowski sum of linear segments
$\{s_i\}_{i=1}^m$:
$$Z=\bigoplus_{i=1}^m \  s_i.$$

For each subset $I\subset [m]$  with $|I|\leq n$, which gives
linearly independent $\{s_i\}_{i\in I}$,  denote by $Z_I$ the
semiopen brick spanned by segments $\{s_i\}_{i\in I}$. It is
well-known that $Z$ can be partitioned into the union of all such
$Z_I$, which implies immediately
$$\Lambda(Z)=\sum_{I\subset [m]}\sharp(Z_I),$$

where $\sharp(\cdot)$ denotes the number of integer points in a
semiopen brick provided that the brick is spanned by linearly
independent vectors. For  linearly dependent vectors we set
$\sharp:=0.$

 For positive integer numbers $\lambda_1,...,\lambda_n$ let us sum up the dilated segments $\lambda_is_i$. Clearly, we have

 $$\Lambda\Big(\bigoplus_{i=1}^m\ \lambda_i s_i\Big)=\sum_{I\subset [m]}\  \sharp(Z_I)\cdot\prod_{i\in I}\lambda_i.$$

 For fixed $s_i$, $\Lambda$ is a polynomial in
 $\lambda_i$, which counts  not only the number of integer points in a convex zonotope (which originates from positive  $\lambda_i$),
 but also the number of integer points in a virtual zonotope, (which originates from any integer $\lambda_i$, including negative ones), see \cite{pkh,panstr}.

\bigskip

 \textbf{Remark.} According to Khovanskii's and Pukhlikov's
 construction \cite{pkh}, given a lattice virtual polytope, each lattice point
 has a \textit{weight}, which is some (possibly, negative) integer number.  The above defined $\Lambda(\cdot)$
 for virtual zonotopes counts the sum of weights.
This fact
 generalizes the \textit{Erchart's reciprocity law }and
 has many other interpretations, such as Riemann-Roch Theorem for toric
 varieties.

\bigskip

We immediately have:
 \begin{lemma}\label{lemmaIntCountCyclo}Let $E=E_n$ be the set  of edges of the complete graph $K_n$.
For the cyclopermutohedron we have:
   $$\Lambda(\mathcal{CP}_{n+1})= \sum_{(I,M) : \ |I| + |M| \leq n-1} (-1)^{|M|} \cdot \sharp\big(\bigoplus_{ I} q_{ij}+ \bigoplus_{ M} r_k\big)$$

 Here $I$ ranges over  subsets of $E$, whereas $M$ ranges over  subsets of
 $[n]$.\qed
 \end{lemma}

Our next aim is to give a formula for one single summand.

\begin{dfn}\label{DefPartialDecForest}

A \textit{partial decorated forest} $F=(G,M)$  is a graph
\newline $G=([n],I)$ without cycles on $n$ labeled vertices together with
a set of marked vertices $M\subset[n]$ such that the following
conditions hold:\begin{enumerate}
    \item Number of marked vertices $|M|$  $+$ number of edges  $|I|$
    is smaller or equal than $n$.
    \item Each connected component of $G$ has at most one marked
    vertex.
\end{enumerate}
\end{dfn}

We already know that decorated forests are in a bijection with
linearly independent $(n-1)$-tuples of $\{q_{ij},r_k\}$ (see Section
\ref{SecVolIsZero}). Therefore, partial decorated forests are in a
bijection with linearly independent collections of segments
$\{q_{ij},r_k\}$.

From now on, we fix one particular  partial decorated forest $F$ and
work with the associated segments.

\bigskip

\textbf{Notation}: $F$ splits into a disjoint union of two forests:
(1) a forest $T=T(F)$ without marked vertices, which is called
\textit{the free forest}, and (2) a rooted forest $R(F)$. In turn,
$T$ is a disjoint union of trees $T_j(F)$.

As in the previous sections,  $C(\cdot)$ denotes the number of
connected components of a forest. In particular, $C(R(F))=|M|$ is
the number of marked vertices.

\bigskip

In this notation we have:

\begin{lemma}\label{LemmaSummand} For the number $\sharp$ of integer points in the semiopen
brick spanned by $\{q_{ij}, r_k\}_{}$, we have:
\begin{enumerate}
    \item If the segments in question do not come from a partial decorated
    forest, then $$\sharp=0.$$
    \item If the segments in question  come from a decorated partial
    forest $F$ with at least one marked vertex, then $$\sharp=n^{|M| - 1} \cdot \gcd[V({T}_{1}),\dots , V({T}_{C(T)} )],$$

    where $|M|$ is the number marked vertices in $F$,  ${T}_{i}$ are the connected
    components of the free forest  $T$, $V({T}_{i})$ is the
    number of vertices in ${T}_{i}$.

    \item If the segments in question  come from a decorated partial
    forest $F$ with no marked vertices, then $$\sharp=1.$$
\end{enumerate}

\end{lemma}

For the proof of the lemma, see Section \ref{LemmaSummand}.\qed

\bigskip

Basing on the lemma, we obtain:

\begin{theorem}\label{ThmInteger}

Define
$$\Phi(v) = \sum_{T} \gcd[\{V({T}_{i})\}],$$
 where the sum ranges over all (non-rooted) forests $T$
on $v$ labeled vertices, $T_i$ are the trees in the forest $T$, and
$V(\cdot)$ is the number of vertices.

\bigskip

Then

$$\Lambda( CP_{n+1} ) =
 \varphi(n) - \sum_{v = 1}^{n-1} \binom{n}{v}  (-v)^{n - v -1}\cdot \Phi(v)= $$

$$=\Lambda(\Pi_n) - \sum_{v = 1}^{n-1} \binom{n}{v}  (-v)^{n - v -1}\cdot \Phi(v), $$

 where $\varphi(n)$ is the number of (non-rooted) forests on $n$
labeled vertices.

\end{theorem}
Proof.
\begin{enumerate}
    \item We count partial decorated forests with no marked points separately. Altogether they
    contribute   $ \varphi(n)=\Lambda(\Pi_n)$.
    \item Next we choose $v$ vertices of the free forest. This can be done in $\binom{n}{v}$ ways.
    \item Each of the forests gives us its own $gcd$. Altogether they give us $
    \Phi(v)$.

    \item Next, we count rooted forests on the remaining  $n-v$ vertices. Each forest is counted with multiplicity  $
    (-n)^{C(F)}$.
 The equality
$$\sum_{f \ is \ a \ rooted \ forest\   \ on \ m \
vertices}x^{C(f)}=x\cdot(x+m)^{m-1}$$ (see Section \ref{SecAbel})
completes the proof. \qed
\end{enumerate}

\textbf{Examples:}

$\Lambda(\mathcal{CP}_{3}) = 1$,

$\Lambda(\mathcal{CP}_{4}) = 18$.

\subsection{Proof of Lemma \ref{LemmaSummand}}

%Technical notation:
We fix a partial decorated forest and the corresponding semiopen
brick spanned by $\{q_{ij}, r_k\}$. The vectors $r_{l}$ will be
called \textit{ long vectors}, whereas $q_{ij}$ will be called
\textit{ short vectors}.
%\textit{The leading
%coordinate} of a long vector $r_{l}$ is the $l^{th}$ coordinate.

%A long vector $r_{i}$ \textit{corresponds to a tree} $P_{j} $ if the
%leading coordinate of $r_{i}$ lies in the vertex  set $Vert(P_{j})$.

\bigskip

As the main tool, we shall use the following  lemma, whose proof
comes from elementary linear algebra.
\begin{lemma}
\label{prop:pm}\begin{enumerate}
    \item The number of integer points
$\sharp(\{v_{i}\}) $ doesn't change if we replace any $v_j$ by the
vector $$v_j + \sum_{i\neq j} (\pm v_i).$$
    \item For an integer $\lambda$, we have:
 $$\sharp(\{\lambda\cdot v_1,v_2,v_3...,v_k\}) =\lambda\cdot\sharp(\{
 v_1,v_2,v_3...,v_k\}).$$
    \item Suppose there exists a coordinate $x_{j}$ such
that among vectors $\{v_{i}\}$ only one vector (say, $v_1$) has
nonzero $j^{th}$ coordinate which  equals $\pm 1$. This will be
called the free coordinate. Then  we can remove $v_1$ from the
collection of segments without changing  the number of integer
points:$$\sharp(\{v_{i}\})=\sharp(\{v_{i}\}_{i \neq 1}).$$
    \item Given a partial decorated forest, replace all the
trees by path trees, keeping for each tree the set of its vertices.
This manipulation does not change the value $\sharp(F)$.
    \item Given one vector $v=(V_1,...,V_n)$,
$$\sharp(v)=\gcd[\{V_i\}],$$
where $\gcd$ denotes the greatest common divisor.\qed
\end{enumerate}
\end{lemma}

\bigskip

\begin{figure}[h]
\centering
\includegraphics[width=12 cm]{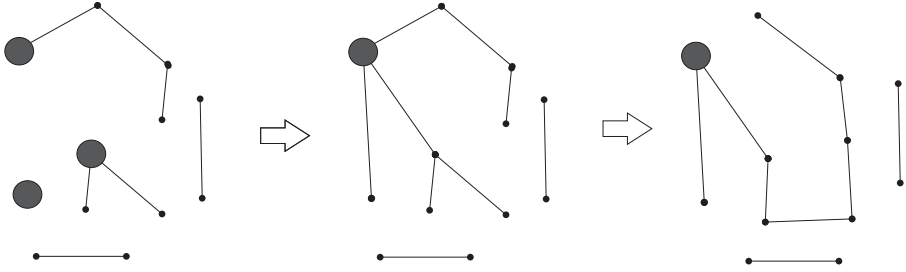}
\caption{Reduction of a partial decorated forests. Grey balls denote
the marked vertices.}\label{RedPart}
\end{figure}

\textit{Reduction of a partial decorated forest} (see Figure
\ref{RedPart}) goes as follows:  Assume we have a partial decorated
forest $F$. \begin{enumerate}
    \item  Choose a marked vertex. We shall call it the \textit{principal marked
vertex}.
    \item Join
the  principal marked vertex   with each of the
other marked vertices by an edge.
    \item  Remove all marks from the
marked vertices that are not principal.
\item Replace the tree with the
marked vertex by a path tree on the same vertices.  We arrive at a
partial decorated forest $\overline{F}$.
\end{enumerate}
Lemma~\ref{prop:pm} implies:
\begin{lemma}
\label{prop:reduction}For a partial decorated forest $F$ and its
reduction $\overline{F}$, we have:

$$\sharp(F)=n^{|M|-1}\cdot \sharp(\overline{F}),$$
where $|M|=|M(F)|$ is the number of marked vertices in $F$.\qed
\end{lemma}

Now we are ready to calculate one single summand from Lemma
\ref{lemmaIntCountCyclo}.  We arrange the column vectors in a
matrix: first come all the $q_{ij}$, after them come all the $r_k$. The main idea is
that the reduction process encodes the way of manipulating with the
columns in the matrix. Using Lemma~\ref{prop:pm}, (1), we can assume
that all the $P_i$ are path trees.
\begin{enumerate}
    \item Assume that the collection contains  some long vector.

The algorithm runs as follows: first, we take the long vector which
corresponds to the principal marked vertex and subtract it from all
the other long vectors. Each of the long vectors (except for the
first one) yields a multiple $n$ and a new short vector.

Next, we subtract the short vectors from the (unique that survived)
long vector aiming at killing its coordinates. Finally, we get a
matrix which allows to remove vectors using Lemma \ref{prop:pm},
(3). Eventually we arrive at
              $$ n^{|M| - 1}\cdot\sharp\begin{pmatrix}
-V(T_{1}) \\
\vdots \\
-V(T_{C(T)})\\
V(T_{1}) +\dots +V(T_{C(T)}) \\
\end{pmatrix}=  n^{|M| - 1}\cdot \gcd[V(T_{1}),\dots, V(T_{C(T)})]. $$

    \item If there are no long vectors in the collection, we remove
    the vectors one by one using Lemma \ref{prop:pm}, (3), and arrive
    at $\sharp=1$.\qed
\end{enumerate}

\begin{figure}[h]
\centering
\includegraphics[width=12 cm]{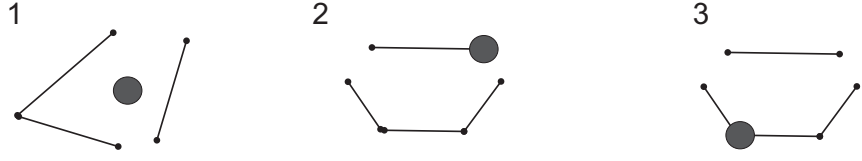}
\caption{Partial decorated forests. Grey balls denote the marked
vertices.}\label{PartialDeco}
\end{figure}

\textbf{Examples.} We exemplify below the reduction for three
collections of vectors. Corresponding partial decorated forests are
depicted in Fig. \ref{PartialDeco}.
\begin{enumerate}
\item  Two free trees with $V=2$ and $V=3$, $|M|=1$.
 $$\sharp\left(
              \begin{array}{cccc}
                1 & 0 & 0 & -1 \\
                -1 & 1 & 0 & -1 \\
                0 & -1 & 0 & -1 \\
                0 & 0 & 1 & -1 \\
                0 & 0 & -1 &-1 \\
                0 & 0 & 0 & 5 \\
              \end{array}
            \right)=\sharp\left(
              \begin{array}{cccc}
                1 & 0 & 0 & 0 \\
                -1 & 1 & 0 & 0 \\
                0 & -1 & 0 & -3 \\
                0 & 0 & 1 & -1 \\
                0 & 0 & -1 &-1 \\
                0 & 0 & 0 & 5 \\
              \end{array}
            \right)=$$

$$=5\cdot\sharp\left(
              \begin{array}{cccc}
                1 & 0 & 0 & 0 \\
                -1 & 1 & 0 & 0 \\
                0 & -1 & 0 & -3 \\
                0 & 0 & 1 & 0 \\
                0 & 0 & -1 &-2 \\
                0 & 0 & 0 & 5 \\
              \end{array}
            \right)=  \sharp\left(
              \begin{array}{c}
                0   \\
                0 \\
               -3    \\
                0    \\
                -2 \\
  5    \\
              \end{array}
            \right) =1.$$

\item One free tree with $V=4$, $|M|=1$.

$$\sharp\left(
              \begin{array}{ccccc}
                1 & 0 & 0&0 &-1 \\
                -1 & 0 & 0 &0& 5 \\
                0 & -1 & 0 &0& -1 \\
                0 & 1 & -1 &0& -1 \\
                0 & 0 & 1 &1& -1 \\
                0 & 0 & 0 &-1& -1 \\
              \end{array}
            \right)= \sharp\left(
              \begin{array}{ccccc}
                1 & 0 & 0&0 &-1 \\
                -1 & 0 & 0 &0& 5 \\
                0 & -1 & 0 &0& -1 \\
                0 & 1 & -1 &0& -1 \\
                0 & 0 & 1 &1& -2 \\
                0 & 0 & 0 &-1& 0 \\
              \end{array}
            \right)= $$

            $$=\sharp\left(
              \begin{array}{ccccc}
                1 & 0 & 0&0 &-1 \\
                -1 & 0 & 0 &0& 5 \\
                0 & -1 & 0 &0& -1 \\
                0 & 1 & -1 &0& -3 \\
                0 & 0 & 1 &1& 0 \\
                0 & 0 & 0 &-1& 0 \\
              \end{array} \right) = \sharp\left(
              \begin{array}{ccccc}
                1 & 0 & 0&0 &-1 \\
                -1 & 0 & 0 &0& 5 \\
                0 & -1 & 0 &0& -4 \\
                0 & 1 & -1 &0& 0 \\
                0 & 0 & 1 &1& 0 \\
                0 & 0 & 0 &-1& 0 \\
              \end{array}
            \right) =$$ $$
           = \sharp\left(
              \begin{array}{ccccc}
                1 & 0 & 0&0 &0 \\
                -1 & 0 & 0 &0& 4 \\
                0 & -1 & 0 &0& -4 \\
                0 & 1 & -1 &0& 0 \\
                0 & 0 & 1 &1& 0 \\
                0 & 0 & 0 &-1& 0 \\
              \end{array}
            \right)= \sharp\left(
              \begin{array}{c}
               0 \\
                4 \\
                -4 \\
                0 \\
                0 \\
                 0 \\
              \end{array}  \right)=4.$$

\item One free tree with $V=2$, $|M|=1$.

$$\sharp\left(
              \begin{array}{ccccc}
                1 & 0 & 0&0 &-1 \\
                -1 & 0 & 0 &0& -1 \\
                0 & -1 & 0 &0& -1 \\
                0 & 1 & -1 &0& 5 \\
                0 & 0 & 1 &1& -1 \\
                0 & 0 & 0 &-1& -1 \\
              \end{array}
            \right)= \sharp\left(
              \begin{array}{ccccc}
                1 & 0 & 0&0 &0 \\
                -1 & 0 & 0 &0& -2 \\
                0 & -1 & 0 &0& -1 \\
                0 & 1 & -1 &0& 5 \\
                0 & 0 & 1 &1& -2 \\
                0 & 0 & 0 &-1& 0 \\
              \end{array}
            \right)= $$

$$ = \sharp\left(
              \begin{array}{cccc}
                1 & 0 & 0&0 \\
                -1 & 0 & 0 & -2 \\
                0 & -1 & 0 & -1 \\
                0 & 1 & -1 & 3 \\
                0 & 0 & 1 & 0 \\
                0 & 0 & 0 & 0 \\
              \end{array}
            \right)= \sharp\left(
              \begin{array}{ccc}
                 0 & 0 &0 \\
                0 & 0 & -2 \\
                -1 & 0 & 0 \\
                1 & -1 &2 \\
                 0 & 1 & 0 \\
                 0 & 0 & 0 \\
              \end{array}
            \right)= \sharp\left(
              \begin{array}{cc}
                  0 \\
                 -2 \\
                  0 \\
                 2 \\
                   0 \\
                   0 \\
              \end{array}
            \right)=2.$$

\end{enumerate}

\end{document}